# On a Laplacian which acts on symmetric tensors

Mikeš J.

Department of Algebra and Geometry, Palacky University,
77146 Olomouc, Listopadu 12, Czech Republic
e-mal: *josef.mikes@upol.cz*

Stepanov S. E., Tsyganok I.I.

Department of Mathematics, Finance University,
125468 Moscow, Leningradsky Prospect, 49-55, Russian Federation
e-mail*: s.e.stepanov@mail.ru, i.i.tsyganok@mail.ru*

**Abstract.** In the present paper we show properties of a little-known Laplacian operator acting on symmetric tensors. This operator is an analogue of the well known Hodge-de Rham Laplacian which acts on exterior differential forms. Moreover, this operator admits the Weitzenböck decomposition and we study it using the analytical method, due to Bochner, of proving vanishing theorems for the null space of a Laplace operator admitting a Weitzenböck decomposition and further of estimating its lowest eigenvalue.

## 1. Introduction

In present paper we show properties of a little-known differential operator $\Delta_{sym}$ of second order which acts on symmetric tensors on a Riemannian manifold $(M, g)$. This operator is an analogue of the well known Hodge-de Rham Laplacian $\Delta$ which acts on exterior differential forms (see [1, p. 54]; [2, p. 204]). The operator $\Delta_{sym}$ is a self-adjoint Laplacian operator and its kernel is a finite-dimensional vector space on a compact Riemannian manifold $(M, g)$. In addition, the Laplacian $\Delta_{sym}$ admits the Weitzenböck decomposition and so we can study it using the analytical method, due to Bochner, of proving vanishing theorems for the null space of a Laplace operator admitting a Weitzenböck decomposition and further of estimating its lowest eigenvalue (see [1, p. 53]; [2, p. 211]; [3]; [4] and [5]).

The paper is organized as follows. In the next section, we give a brief review of the Riemannian geometry of the Laplacian $\Delta_{sym}$. In the third section of the paper we prove vanishing theorems of the ker-nel of $\Delta_{sym}$. In the fourth section, we demonstrate the main result of the present paper which is the estimation of the first eigenvalue of the Laplacian $\Delta_{sym}$. In the final section a Bochner integral formula for the Yano Laplacian will be obtained.

## 2. Definitions and notations

**2.1**. Let $(M, g)$ be a compact oriented $C^\infty$-Riemannian manifold of dimension $n \geq 2$ and $S^p M$ be a symmetric tensor product of order $p \geq 1$ of a tangent bundle $TM$ of $M$. On tensor spaces on $M$ we have *canonical scalar product* $g(\,\cdot\,,\,\cdot\,)$ and on their $C^\infty$-sections the *global scalar product* $\langle\,\cdot\,,\,\cdot\,\rangle$. In particular, for any $\varphi, \varphi' \in C^\infty S^p M$ we have

$$\langle \varphi, \psi \rangle = \int_M \frac{1}{p!} g(\varphi, \varphi') \, dv, \tag{2.1}$$

where $dv$ is the volume element of $(M, g)$. If $g_{ij}, \varphi_{i_1 \ldots i_p}$ and $\varphi'_{i_1 \ldots i_p}$ denote components of the metric $g$ and tensor fields $\varphi, \varphi' \in C^\infty S^p M$ with respect to a local coordinate system $x^1, \ldots, x^n$ on $(M, g)$ then

$$\langle \varphi, \psi \rangle = \int_M \frac{1}{p!} \left( \varphi^{j_1 \ldots j_p} \varphi'_{j_1 \ldots j_p} \right) dv$$

where $\varphi^{j_1 \ldots j_p} = g^{i_1 j_1} \ldots g^{i_p j_p} \varphi_{i_1 \ldots i_p}$ and $g^{ik} g_{jk} = \delta^i_j$ for the Kronecker delta $\delta^i_j$.

Next, if $D$ is a differential operator between some tensor bundles over $M$, its formal adjoint $D^*$ is uniquely defined by the formula $\langle D\,\cdot\,,\,\cdot\,\rangle = \langle\,\cdot\,, D^*\,\cdot\,\rangle$ (see [1, p. 460]). For example (see [1, p. 514]), the covariant derivative $\nabla : C^\infty S^p M \to C^\infty(T^* M \otimes S^p M)$, whose formal adjoint $\nabla^*$ will be also denoted by $\delta = \nabla^* : C^\infty(T^* M \otimes S^p M) \to C^\infty S^p M$. In local coordinates, this operator is defined by equalities $(\delta \xi)_{i_1 \ldots i_p} = -g^{ij} \nabla_j \xi_{i i_1 \ldots i_p}$ for an arbitrary $\xi \in C^\infty(T^* M \otimes S^p M)$.

Furthermore we define (see also [1, p. 514]) the operator $\delta^* : C^\infty S^p M \to C^\infty S^{p+1} M$ of degree 1 by the formula $\delta^* \varphi = \mathrm{Sym}(\nabla \varphi)$ for an arbitrary $\varphi \in C^\infty S^p M$. In local coordinates we have

$$(\delta^* \varphi)_{i_0 i_1 \ldots i_{p-1} i_p} = \nabla_{i_0} \varphi_{i_1 \ldots i_{p-1} i_p} + \ldots + \nabla_{i_p} \varphi_{i_0 i_1 \ldots i_{p-1}}.$$

Then $\delta: C^\infty S^{p+1}M \to C^\infty S^p M$ is the adjoint operator of $\delta^*$ with respect to the global product (2.1).

**2.2**. In [6] the operator $\Delta_{sym} = \delta\delta^* - \delta^*\delta : C^\infty S^p M \to C^\infty S^p M$ of degree 2 was defined. The operator $\Delta_{sym}$ is related to a variational problem (see [6]), as follows: If we define the "energy" of symmetric tensor field $\varphi$ by $E(\varphi) = \frac{1}{2}\langle \varphi, \Delta_{sym}\varphi \rangle$, then $\Delta_{sym}\varphi = 0$ is the condition for a free extremal of $E(\varphi) = \frac{1}{2}\langle \varphi, \Delta_{sym}\varphi \rangle$. This operator is an analogue of the well known Hodge-de Rham Laplacian $\Delta: C^\infty \Lambda^p M \to C^\infty \Lambda^p M$ which acts on $C^\infty$-sections of the bundle $\Lambda^p M$ of covariant skew-symmetric tensors of degree $p$; in other words, exterior differential $p$-forms on $M$ (see [1, p. 34]; [2, p. 204]).

The operator $\Delta_{sym}$ is a self-adjoint Laplacian operator and its kernel is a finite-dimensional vector space on a compact Riemannian manifold $(M, g)$ (see [6]; [7], [10]). In this case, by virtue of the *Fredholm alternative* (see [1, p. 464]; [2, p. 205]) vector spaces Ker $\Delta_{sym}$ and Im $\Delta_{sym}$ are orthogonal complements of each other with respect to the global scalar product (2.1), i.e.

$$C^\infty S^p M = \text{Ker } \Delta_{sym} \oplus \text{Im } \Delta_{sym}. \tag{2.2}$$

At the same time, we recall that this decomposition is an analog of the well-known *Hodge orthogonal decomposition* $C^\infty \Lambda^p M = \text{Ker } \Delta \oplus \text{Im } \Delta$ where $1 \leq p \leq n-1$ and the kernel of $\Delta$ consists of *harmonic* exterior differential $p$-forms on $M$ (see [2], pp. 202-213).

Compare the operator $\Delta_{sym}$ with the *rough Laplacian* $\nabla^*\nabla$ (see [1, p. 54]). First, it is easy to see that these two operators coincide if $(M, g)$ is a locally Euclidean space. Second, the operator $\Delta_{sym} - \nabla^*\nabla$ has order zero and can be defined by symmetric endomorphism $B_p$ of the bundle $S^p M$ where $B_p$ can be algebraically (even linearly) expressed through the curvature and Ricci tensors of $(M, g)$ (see [1, p. 53]; [6]). This can be expressed by the Weitzenböck decomposition formula $\Delta_{sym} = \nabla^*\nabla - B_p$ and means that $\Delta_{sym}$ is a *Lichnerowicz Laplacian*. In local coordinates we have

$$B_p(\varphi)_{i_1 \ldots i_p} = \sum_k g^{js} r_{i_k j} \varphi_{i_1 \ldots i_{k-1} s i_{k+1} \ldots i_p} - \sum_{k \neq l} g^{js} g^{mt} R_{i_k j i_l m} \varphi_{i_1 \ldots i_{k-1} s i_{k+1} \ldots i_{l-1} t i_{l+1} \ldots i_p}, \tag{2.3}$$

where $r_{ij}$ and $R_{ijkl}$ are components of the Ricci tensor *Ric* and curvature tensor *R* with respect to a local coordinate system $x^1,...,x^n$.

**Remark**. Eigenvalues with their multiplicities of the Lichnerowicz Laplacian $\Delta_L = \nabla^*\nabla + B_p$, which acts on the vector space $C^\infty S^p M$ over the Euclidian sphere $S^n$ were discussed in [11].

In particular, if (*M, g*) is a locally Euclidean manifold then the equation $\Delta_{sym}\varphi = 0$ becomes $\sum_k \frac{\partial^2 \varphi_{i_1...i_p}}{(\partial x^k)^2} = 0$ with respect to a local Cartesian coordinate system $x^1,...,x^n$. This means that all components of a covariant completely symmetric *p*-tensor $\varphi$ on a locally Euclidean manifold (*M, g*) are harmonic functions belonging to the kernel of the Laplacian $\Delta_{sym}$. This property is also characteristic for any harmonic exterior differential *p*-form (see [2, p. 205-212]). Therefore, the symmetric tensor $\varphi \in \ker \Delta_{sym}$ was named in [6] as a *harmonic symmetric p-tensor* on a Riemannian manifold (*M, g*).

**Remark.** For a harmonic symmetric *p*-tensor, it was said in [6] that it does not have any geometrical interpretation for $p \geq 2$. It was also pointed out that symmetric tensors, whose covariant derivative vanishes is clearly harmonic, the metric tensor *g* being the most important example. These tensors are trivial examples of a harmonic symmetric tensor. Now we can say (see Theorem 4.6 in the paper [12] of the first of the three authors of this paper) that the 3-tensor *T* of a nearly integrable *SO*(3)-structure on a five-dimensional connected Riemannian manifold (*M, g*) is a non-trivial example of a harmonic symmetric tensor.

**2.3**. For $p \geq 2$, the explicit expression for $\Delta_{sym}$ is sufficiently complicated but for the case $p = 1$, it has the form $\Delta_{sym} = \nabla^*\nabla - B_1$ where $B_1 = Ric$. From well known Weitzenböck decomposition formula $\Delta = \nabla^*\nabla + Ric$ we obtain $\Delta_{sym} = \Delta - 2Ric$, where, in local coordinates, we have $(\Delta_{sym}\varphi)_i = (\Delta\varphi)_i - 2r_{ij}\varphi^j$ (see [7]; [10]). This form of the operator $\Delta_{sym}$ was used by K. Yano (see, for example, the monograph [8, p. 40] and paper [9]) for the investigation of *infinitesimal isometries* of

($M, g$) preserving the metric $g$ and special vector fields named *geodesic vector fields*. Therefore, in paper [10], we called it as the *Yano operator* or, in other words, the *Yano Laplacian*. We will use this name for $\Delta_{sym}$ throughout this paper.

On the other hand, a vector field $X$ on ($M, g$) is called an *infinitesimal harmonic transformation* if the one-parameter group of local transformations of ($M, g$) generated by $X$ consists of local harmonic diffeomorphisms (see [13]). At the same time, the equality $\Delta_{sym} \varphi = 0$ for $\varphi = g(X, \cdot)$ is necessary and sufficient condition for a vector field $X$ to be an infinitesimal harmonic transformation on ($M, g$) (see [7], [10]). Hence, in the special case of $p = 1$ the kernel of $\Delta_{sym}$ consists of infinitesimal harmonic transformations.

In addition, we recall the following examples of infinitesimal harmonic transformations: infinitesimal isometric transformations of Riemannian manifolds (see [8, p. 44]), infinitesimal conformal transformations of two-dimensional Riemannian manifolds (see [8, p. 47]), affine Killing (see [8, p. 45]) and geodesic vector fields (see [9]) on Riemannian manifolds, holomorphic vector fields on nearly Kählerian manifolds (see [10]) and vector fields that transform a Riemannian metrics into *Ricci soliton* metrics (see [7]).

In conclusion, we point out that the vector space $\operatorname{Ker} \delta^* := \operatorname{Ker} \Delta_{sym} \cap \operatorname{Ker} \delta$ of all infinitesimal isometric transformations must be an orthogonal complement of the vector space $\operatorname{Ker} \Delta_{sym} \cap \operatorname{Im} d$ of all infinitesimal gradient harmonic transformations with respect to the whole space $\operatorname{Ker} \Delta_{sym}$ on a compact orientable Einstein manifold ($M, g$) (see [9]). In view of (2.2) we can write the following orthogonal decomposition $C^\infty S^1 M = \operatorname{Ker} \Delta_{sym} \cap \operatorname{Im} d \oplus \operatorname{Ker} \delta^* \oplus \operatorname{Im} \Delta_{sym}$.

### 3. Vanishing theorems

**3.1**. Let ($M, g$) be a compact and oriented Riemannian manifold of dimension $n \geq 2$ and $\Delta_{sym} : C^\infty S^p M \to C^\infty S^p M$ be the Yano Laplacian. Then the following proposition is true (for the proof see Appendix).

**Lemma 3.1**. *Let (M, g) be a compact oriented $C^\infty$-Riemannian manifold of dimension $n \geq 2$ and $\Delta_{sym} : C^\infty S^p M \to C^\infty S^p M$ be the Yano Laplacian then*

$$\langle \Delta_{sym} \varphi, \varphi \rangle = \langle \nabla \varphi, \nabla \varphi \rangle - \langle B_p(\varphi), \varphi \rangle \qquad (3.1)$$

where $\langle B_p(\varphi), \varphi \rangle = \int_M \frac{1}{(p-1)!} g(B_p(\varphi), \varphi) dv$.

The integral formula (3.1) we obtain by the Weitzenböck decomposition formula $\Delta_{sym} = \nabla^* \nabla - B_p$. In particular, (2.3) becomes

$$\langle B_p(\varphi), \varphi \rangle = \int_M \frac{1}{(p-1)!} \left( r_{ij} \varphi^{i i_2 \ldots i_p} \varphi^{j}{}_{i_2 \ldots i_p} - (p-1) R_{ijkl} \varphi^{i k i_3 \ldots i_p} \varphi^{jl}{}_{i_3 \ldots i_p} \right) dv$$

with respect to a local coordinate system $x^1, \ldots, x^n$.

For the special case of $p = 1$ we can rewrite the formula (3.1) in the form

$$\langle \Delta_{sym} \varphi, \varphi \rangle = \langle \nabla \varphi, \nabla \varphi \rangle - \langle Ric\, X, X \rangle \qquad (3.2)$$

where $\varphi = g(X, \cdot)$ for an arbitrary $X \in C^\infty T^* M$. It is well known that $\Delta = \nabla^* \nabla + Ric$. Then from (3.2) we obtain the formula

$$\langle \Delta_{sym} \varphi, \varphi \rangle = \langle \Delta \varphi, \varphi \rangle - 2 \langle Ric\, X, X \rangle . \qquad (3.3)$$

**Corollary 3.2.** *Let (M, g) be a compact oriented $C^\infty$-Riemannian manifold of dimension $n \geq 2$ and $\Delta_{sym} : C^\infty S^1 M \to C^\infty S^1 M$ be the Yano Laplacian then $\langle \Delta \varphi, \varphi \rangle - \langle \Delta_{sym} \varphi, \varphi \rangle = 2 \langle Ric\, X, X \rangle$ where $\varphi = g(X, \cdot)$ for an arbitrary $X \in C^\infty T^* M$.*

**3.2.** The curvature tensor $R$ defines a linear endomorphism $\overset{\circ}{R} : S^2(M) \to S^2(M)$ by the formula $\overset{\circ}{R}(\varphi)_{il} = R_{ijkl} \varphi^{jk}$ where $R_{ijkl}$ are local components of $R$ and $\varphi^{ij} = g^{ik} g^{jl} \varphi_{kl}$ for local components $\varphi_{ij}$ of an arbitrary $\varphi \in C^\infty S^2 M$. The symmetries of $R$ imply that $\overset{\circ}{R}$ is a self-adjoint operator, with respect to the pointwise inner product on $S^2(M)$. Hence the eigenvalues of $\overset{\circ}{R}$ are all real numbers at each point $x \in M$. Thus, we say $\overset{\circ}{R}$ is positive (resp. negative), or simply $\overset{\circ}{R} > 0$ (resp. $\overset{\circ}{R} < 0$), if all the eigenvalues of $\overset{\circ}{R}$ are positive (resp. negative). This endomorphism $\overset{\circ}{R}: S^2(M) \to S^2(M)$ was named as the *operator of curvature of the second kind* (see [15]).

**Remark.** Definition, properties and applications of $\overset{\circ}{R}$ can be found in the following monograph and papers [1]; [14]; [15]; [16]; [17] and etc.

Let $(M, g)$ be an $n$-dimensional Riemannian manifold. If there exists a positive constant $\varepsilon$ such that

$$R_{ijkl}\, \varphi^{il}\, \varphi^{jk} \leq -\varepsilon\, \varphi_{ij}\, \varphi^{ij} \tag{3.4}$$

holds for any $\varphi \in C^\infty S^2 M$, then $(M, g)$ is said to be a *Riemannian manifold of negative restricted curvature operator of the second kind*. In this case the following theorem is true.

**Theorem 3.3.** *Let $(M, g)$ be an n-dimensional compact oriented Riemannian manifold. Suppose the curvature operator of the second kind $\overset{\circ}{R}: S^2(M) \to S^2(M)$ is negative and bounded above by some negative number $-\varepsilon$ at each point $x \in M$. Then $(M, g)$ does not admit non-zero symmetric harmonic p-tensors.*

**Proof.** Let $\varphi$ be a harmonic symmetric $p$-tensor then from (3.1) we obtain

$$\langle B_p(\varphi), \varphi \rangle = \langle \nabla \varphi, \nabla \varphi \rangle \geq 0. \tag{3.5}$$

On the other hand, from (3.4) we obtain the following inequalities

$$-(p-1)R_{ijkl}\, \varphi^{iki_3...i_p}\, \varphi^{jl}{}_{i_3...i_p} = (p-1)R_{ijlk}\, \varphi^{iki_3...i_p}\, \varphi^{jl}{}_{i_3...i_p} \leq -(p-1)\varepsilon\, \varphi_{i_1...i_p}\, \varphi^{i_1...i_p};$$

$$r_{ij}\, \varphi^{ii_2...i_p}\, \varphi^{j}{}_{i_2...i_p} \leq -(n-1)\varepsilon\, \varphi_{i_1...i_p}\, \varphi^{i_1...i_p},$$

which show that

$$\langle B_p(\varphi), \varphi \rangle \leq -p\,(n+p-2)\langle \varphi, \varphi \rangle < 0 \tag{3.6}$$

for an arbitrary non-zero $\varphi \in C^\infty S^p M$. Inequalities (3.5) and (3.6) contradict each other. This contradiction proves our theorem.

**Remark.** We recall that an orientable and compact Riemannian manifold (M, g) does not admit a non-zero harmonic exterior differential p-form if its curvature operator $\overset{\circ}{R}: S^2(M) \to S^2(M)$ is positive and bounded below (see [15]).

**3.3.** A 1-dimensional immersed submanifold $\gamma$ of (M, g) is called a *geodesic* if exists a parameterization $\gamma: x^i = x^i(t)$ for $t \in I \subset \mathbb{R}$ satisfying $\nabla_{\dot{x}} \dot{x} = 0$. If each solution $x^i = x^i(t)$ of the equations $\nabla_{\dot{x}} \dot{x} = 0$ of the geodesics satisfies the condition $\varphi(\dot{x},...,\dot{x}) = const$ for smooth covariant completely symmetric p-tensor $\varphi$ and $\dot{x} = \frac{d x^k}{dt} \cdot \frac{\partial}{\partial x^k}$, then the equations $\nabla_{\dot{x}} \dot{x} = 0$ admit so called a *first integral of the p-th order of differential equations of geodesics*.

The equation $(\delta^* \varphi)(X, X, ..., X) = (\nabla_X \varphi)(X, X, ..., X) = 0$ serves as necessary and sufficient condition for this [18, pp. 128-129]. On the other hand, the integral formula (3.1) can be given in the form

$$\langle B_p(\varphi), \varphi \rangle + \langle \delta^* \varphi, \delta^* \varphi \rangle - \langle \delta \varphi, \delta \varphi \rangle - \langle \nabla \varphi, \nabla \varphi \rangle = 0. \tag{3.7}$$

Thus, by the Theorem 3.1, we have the following result (see also [19]).

**Corollary 3.4.** *Let (M, g) be an n-dimensional compact oriented Riemannian manifold. Suppose the curvature operator of the second kind $\overset{\circ}{R}: S^2(M) \to S^2(M)$ is negative and bounded above by some negative number $-\varepsilon$ at each point $x \in M$. Then (M, g) does not admit first integrals of the p-th order of differential equations of geodesics.*

Finally, for the special case of p =1, by the Corollary 3.2, we have (see also [10])

**Corollary 3.5.** *If, in n-dimension $(n \geq 2)$ compact orientable Riemannian manifold (M, g), the Ricci tensor Ric is negative definite, then (M, g) does not admit non-zero infinitesimal harmonic transformations.*

# 4. Spectral properties of the Yano Laplacian

**4.1.** A real number $\lambda^r$, for which there is a tensor $\varphi \in C^\infty S^p M$ (not identically zero) such that $\Delta_{sym}\varphi = \lambda^p \varphi$, is called an *eigenvalue* of $\Delta_{sym}$ and the corresponding $\varphi \in C^\infty S^p M$ an *eigentensor* of $\Delta_{sym}$ corresponding to $\lambda^p$. The zero p-tensor and the eigentensors corresponding to a fixed eigenvalue $\lambda^p$ form a vector subspace of $S^p M$ denoted by $V_{\lambda^p}(M)$ and called the eigenspace corresponding to the eigenvalue $\lambda^p$. The following theorem is true.

**Theorem 4.1.** *Let $(M, g)$ be an n-dimensional $(n \geq 2)$ compact and oriented Riemannian manifold and $\Delta_{sym}: C^\infty S^p M \to C^\infty S^p M$ be the Yano Laplacian.*

1) *Suppose the curvature operator of the second kind $\overset{\circ}{R}: S^2(M) \to S^2(M)$ is negative and bounded above by some negative number $-\varepsilon$ at each point $x \in M$. Then an arbitrary eigenvalue $\lambda^p$ of $\Delta_{sym}$ is positive.*

2) *The eigenspaces of $\Delta_{sym}$ are finite dimensional.*

3) *The eigenforms corresponding to distinct eigenvalues are orthogonal.*

**Proof.** 1) Let $\varphi \in V_{\lambda^p}(M)$ be a non-zero eigentensor corresponding to the eigenvalue $\lambda^p$, that is $\Delta_{sym}\varphi = \lambda^p \varphi$ then we can rewrite the formula (3.1) in the form

$$\lambda^p \langle \varphi, \varphi \rangle = -\langle B_p(\varphi), \varphi \rangle + \langle \nabla \varphi, \nabla \varphi \rangle. \tag{4.1}$$

If we suppose that the curvature operator of the second kind $\overset{\circ}{R}: S^2(M) \to S^2(M)$ is negative and bounded above by some negative number $-\varepsilon$ at each point $x \in M$ then, by the inequality (3.1), we conclude

$$\lambda^p \langle \varphi, \varphi \rangle = -\langle B_p(\varphi), \varphi \rangle + \langle \nabla \varphi, \nabla \varphi \rangle \geq p(n+p-2)\langle \varphi, \varphi \rangle + \langle \nabla \varphi, \nabla \varphi \rangle > 0. \tag{4.2}$$

2) The eigenspaces of $\Delta_{sym}$ are finite dimensional because $\Delta_{sym}$ is an elliptic operator.

3) Let $\lambda_1^p \neq \lambda_2^p$ and $\varphi_1, \varphi_2$ be the corresponding symmetric eigentensors. Then $\langle \Delta_{sym} \varphi_1, \varphi_2 \rangle = \lambda_1^p \langle \varphi_1, \varphi_2 \rangle$ and $\langle \Delta_{sym} \varphi_1, \varphi_2 \rangle = \langle \varphi_1, \lambda_2^p \varphi_2 \rangle = \lambda_2^p \langle \varphi_1, \varphi_2 \rangle$. Therefore $0 = (\lambda_1^p - \lambda_2^p)\langle \varphi_1, \varphi_2 \rangle$ and since $\lambda_1^p \neq \lambda_2^p$ it follows that $\langle \varphi_1, \varphi_2 \rangle = 0$, that is, $\varphi_1$ and $\varphi_2$ are orthogonal.

Using the general theory of elliptic operators on compact (M, g) it can be proved that $\Delta_{sym}$ has a discrete spectrum, denoted by $\mathrm{Spec}^{(p)}(M)$, consisting of real eigenvalues of finite multiplicity which accumulate only at infinity. In symbols, we have $\mathrm{Spec}^{(p)}(M) = \{0 \leq |\lambda_1^p| \leq |\lambda_2^p| \leq \ldots \to +\infty\}$. In addition, if we suppose that the curvature operator of the second kind $\overset{\circ}{R}: S^2(M) \to S^2(M)$ is negative and bounded above then $\mathrm{Spec}^{(p)}(M) = \{0 < \lambda_1^p \leq \lambda_2^p \leq \ldots \to +\infty\}$. Moreover, the following theorem is true.

**Theorem 4.2.** *Let (M, g) be an n-dimensional compact and oriented Riemannian manifold. Suppose the curvature operator of the second kind $\overset{\circ}{R}: S^2(M) \to S^2(M)$ is negative and bounded above by some negative number $-\varepsilon$ at each point $x \in M$. Then the first eigenvalue $\lambda_1^p$ of the Yano Laplacian $\Delta_{sym}: C^{\infty} S^p M \to C^{\infty} S^p M$ satisfies the inequality $\lambda_1^p > p(n + p - 2)\varepsilon$.*

**Proof.** Suppose the curvature operator of the second kind $\overset{\circ}{R}: S^2(M) \to S^2(M)$ is negative and bounded above by some negative number $-\varepsilon$ at each point $x \in M$. Then for an eigentensor $\varphi$ corresponding to an eigenvalue $\lambda^p$, (4.2) becomes the inequalities

$$\lambda^p \langle \varphi, \varphi \rangle \geq p(n + p - 2)\varepsilon \langle \varphi, \varphi \rangle + \langle \nabla \varphi, \nabla \varphi \rangle \geq p(n + p - 2)\varepsilon \langle \varphi, \varphi \rangle \qquad (4.3)$$

which prove that

$$\lambda_1^p \geq p(n + p - 2)\varepsilon > 0. \qquad (4.4)$$

If the equality is valid in (4.4), then from (4.3) we obtain $\nabla \varphi = 0$ and $B_p(\varphi) = 0$. In this case we have the identity $\varphi = 0$, because $B_p(\varphi)$ is strictly negative for an arbitrary nonzero $\varphi \in C^{\infty} S^p M$. This completes the proof.

In the special case of $p = 1$ we obtain the following corollary from the above theorem.

**Corollary 4.3**. *Let $(M, g)$ be an n-dimensional compact oriented Riemannian manifold. Suppose the Ricci tensor satisfies the inequality $Ric \leq -(n-1)\varepsilon\, g$ for some positive number $\varepsilon$, then the first eigenvalue $\lambda_1^1$ of the Yano Laplacian $\Delta_{sym}: C^\infty S^1 M \to C^\infty S^1 M$ satisfies the inequality $\lambda_1^1 \geq 2(n-1)\varepsilon$. The equality $\lambda_1^1 = 2(n-1)\varepsilon$ is attained for some harmonic eigenform $\varphi \in C^\infty S^1 M$ and in this case the multiplicity of $\lambda_1^1$ is less than or equals to the Betti number $b_1(M)$.*

**Proof.** Let $(M, g)$ be an $n$-dimensional compact and oriented Riemannian manifold, and define the Ricci tensor *Ric* of $(M, g)$ by its components $r_{ij}$ with respect to a local coordinate system $x^1,...,x^n$. If there exists a positive constant $\varepsilon$ such that

$$r_{ij}\, \varphi^i\, \varphi^j \leq -(n-1)\varepsilon\, \varphi_i\, \varphi^i \tag{4.5}$$

holds for any $\varphi \in C^\infty S^1 M$, then $(M, g)$ is said to be a *Riemannian manifold of negative restricted Ricci curvature*.

Next, suppose that the Ricci is negative and bounded above by some negative number $-(n-1)\varepsilon$ at each point $x \in M$, i.e. $Ric \leq -(n-1)\varepsilon\, g$, then from $\Delta_{sym} = \Delta - 2Ric$ we obtain the inequality

$$\langle \Delta_{sym}\, \varphi, \varphi \rangle \geq 2(n-1)\varepsilon \langle \varphi, \varphi \rangle + \langle \Delta\varphi, \varphi \rangle \tag{4.6}$$

for an arbitrary $\varphi \in C^\infty S^1 M$. Then for an eigenform $\varphi$ corresponding to an eigenvalue $\lambda^1$, (4.6) becomes the inequalities

$$\lambda^1 \langle \varphi, \varphi \rangle \geq 2(n-1)\varepsilon \langle \varphi, \varphi \rangle + \langle \Delta\varphi, \varphi \rangle \geq 2(n-1)\varepsilon \langle \varphi, \varphi \rangle \tag{4.7}$$

which prove that

$$\lambda_1^1 \geq 2(n-1)\varepsilon > 0. \tag{4.8}$$

If the equality is valid in (4.8), then from (4.7) we obtain $\Delta\varphi = 0$. In this case $\varphi$ is a harmonic 1-form and the multiplicity of $\lambda_1^1$ is less than or equals to the Betti number $b_1(M)$.

## 4. Appendix

In this section we prove the Lemma 3.1. For this we introduce a vector field $X$ with components $\varphi^j{}_{i_2...i_p} \nabla_j \varphi^{ii_2...i_p}$ with respect to a local coordinate system $x^1,...,x^n$. Then

$$\operatorname{div} X = \nabla_i \left( \varphi^j{}_{i_2...i_p} \nabla_j \varphi^{ii_2...i_p} \right) = \nabla_i \varphi_{ji_2...i_p} \nabla^j \varphi^{ii_2...i_p} + \varphi^j{}_{i_2...i_p} \nabla_i \nabla_j \varphi^{ii_2...i_p} =$$

$$= \nabla_i \varphi_{ji_2...i_p} \nabla^j \varphi^{ii_2...i_p} + r_{ij} \varphi^{ii_2...i_p} \varphi^j{}_{i_2...i_p} - (p-1) R_{ijkl} \varphi^{iki_3...i_p} \varphi^{jl}{}_{i_3...i_p}.$$

In turn, for the vector field $Y$ with local components $\varphi^i{}_{i_2...i_p} \nabla_j \varphi^{ji_2...i_p}$, we have the equalities

$$\operatorname{div} Y = \nabla_i \left( \varphi^i{}_{i_2...i_p} \nabla_j \varphi^{ji_2...i_p} \right) = \nabla_i \varphi^i{}_{i_2...i_p} \nabla_j \varphi^{ji_2...i_p} + \varphi^i{}_{i_2...i_p} \nabla_i \nabla_j \varphi^{ji_2...i_p}.$$

Let $(M, g)$ be a compact and oriented Riemannian manifold then it follows from the Green's theorem $\int_M \operatorname{div} Z \, dv = 0$ for $Z = X - Y$ that

$$\int_M \left( r_{ij} \varphi^{ii_2...i_p} \varphi^j{}_{i_2...i_p} - (p-1) R_{ijkl} \varphi^{iki_3...i_p} \varphi^{jl}{}_{i_3...i_p} + \right.$$

$$\left. + \nabla_{i_0} \varphi_{i_1 i_2...i_p} \nabla^{i_1} \varphi^{i_0 i_2...i_p} - \nabla_i \varphi^i{}_{i_2...i_p} \nabla_j \varphi^{ji_2...i_p} \right) dv = 0 \qquad (4.1)$$

Since

$$\nabla_{i_0} \varphi_{i_1 i_2...i_p} \nabla^{i_1} \varphi^{i_0 i_2...i_p} = \frac{1}{p(p+1)} \left( \nabla_{i_0} \varphi_{i_1 i_2...i_p} + ... + \nabla_{i_p} \varphi_{i_0 i_1...i_{p-1}} \right) \left( \nabla^{i_0} \varphi^{i_1 i_2...i_p} + ... + \nabla^{i_p} \varphi^{i_0 i_1...i_{p-1}} \right) -$$

$$- \frac{1}{p} \left( \nabla_{i_0} \varphi_{i_1 i_2...i_p} \right) \left( \nabla^{i_0} \varphi^{i_1 i_2...i_p} \right)$$

the integral formula (4.1) can be written in the form

$$\langle B_p(\varphi), \varphi \rangle + \langle \delta^* \varphi, \delta^* \varphi \rangle - \langle \delta \varphi, \delta \varphi \rangle - \langle \nabla \varphi, \nabla \varphi \rangle = 0 \qquad (4.2)$$

where

$$\langle B_p(\varphi), \varphi \rangle = \int_M \frac{1}{(p-1)!} g(B_p(\varphi), \varphi) dv = \int_M \frac{1}{(p-1)!} \left( r_{ij} \varphi^{ii_2...i_p} \varphi^j{}_{i_2...i_p} - (p-1) R_{ijkl} \varphi^{iki_3...i_p} \varphi^{jl}{}_{i_3...i_p} \right) dv;$$

$$\langle \delta^* \varphi, \delta^* \varphi \rangle = \int_M \frac{1}{(p+1)!} \left( \nabla_{i_0} \varphi_{i_1 i_2...i_p} + ... + \nabla_{i_p} \varphi_{i_0 i_1...i_{p-1}} \right) \left( \nabla^{i_0} \varphi^{i_1 i_2...i_p} + ... + \nabla^{i_p} \varphi^{i_0 i_1...i_{p-1}} \right) dv;$$

$$\langle \delta\varphi, \delta\varphi \rangle = \int_M \frac{1}{(p-1)!} \left( \nabla_i \varphi^i{}_{i_2\ldots i_p} \right)\left( \nabla_j \varphi^{ji_2\ldots i_p} \right) dv; \quad \langle \nabla\varphi, \nabla\varphi \rangle = \int_M \frac{1}{p!} \left( \nabla_{i_0} \varphi_{i_1 i_2\ldots i_p} \right)\left( \nabla^{i_0} \varphi^{i_1 i_2\ldots i_p} \right) dv.$$

The operator $\Delta_{sym}$ satisfies the property $\langle \Delta_{sym}\varphi, \varphi \rangle = \langle \delta^*\varphi, \delta^*\varphi \rangle - \langle \delta\varphi, \delta\varphi \rangle$, which follows immediately from its definition. Therefore, the integral formula (4.2) can be rewritten in the form (3.1).